\input amstex
\font\smallcaps=cmcsc10
\documentstyle{amsppt}

\def\bs{\bigskip}

\nologo
\def\un{\underbar}
\def\und{\underline}

\def\l[{\leftbracket}
\def\r]{\rightbracket}
\NoBlackBoxes
\NoRunningHeads

\def\r{\Cal R}
\def\u{\Cal U}

\def\V{\Cal V}
\input epsf

\magnification=\magstep1
\pageheight{7in}

\TagsOnRight
\topmatter
\title {NEW WEIGHTED ROGERS-RAMANUJAN PARTITION THEOREMS 
AND THEIR IMPLICATIONS}\endtitle
\author {\it Krishnaswami Alladi\footnote{Research supported in part by the 
National Science Foundation Grant DMS 0088975} and Alexander 
Berkovich\footnote{Research supported in part by a University of Florida 
CLAS Research Award}}\endauthor
\affil University of Florida, Gainesville, Florida  32611 \endaffil
\abstract
This paper has a two-fold purpose. First, by considering a reformulation 
of a deep theorem of G\"ollnitz, we obtain a new weighted partition 
identity involving the Rogers-Ramanujan partitions, namely, partitions 
into parts differing by at least two.  Consequences of this include 
Jacobi's celebrated triple product identity for theta functions, 
Sylvester's famous refinement of Euler's theorem, as well as certain 
weighted partition identities. Next, by studying partitions with 
prescribed bounds on successive ranks and replacing these with weighted 
Rogers-Ramanujan partitions, we obtain two new sets of theorems - a set 
of three theorems involving partitions into parts 
$\not\equiv 0, \pm i$ ($mod$ 6), and a set of three theorems involving 
partitions into parts  $\not\equiv 0, \pm i$ ($mod$ 7), $i=1,2,3$.
\endabstract
\endtopmatter
\document
\baselineskip=17pt
\head \S0 {\un{Introduction}}\endhead
By a Rogers-Ramanujan partition we mean a partition into parts 
differing by $\ge2$. In this paper we obtain several new results 
by attaching weights to Rogers-Ramanujan partitions. In all instances 
the weights are defined multiplicatively.

In the first part of the paper(\S1 - \S6), we obtain a reformulation 
of a deep theorem of G\"ollnitz [15] and discuss the implications. The 
new reformulation is stated as Theorem 1 in \S1. Special cases of this 
yield Jacobi's triple product identity for theta functions (see \S2), 
Sylvester's famous refinement of Euler's theorem (see \S3), and a 
weighted partition identity connecting partitions into distinct parts 
and Rogers-Ramanujan partitions (see \S4). The proof of Theorem 1 is given 
in \S5-\S6, with \S5 describing the necessary prerequisites, namely, 
the method of weighted words of Alladi-Andrews-Gordon [5], and \S6 giving 
the details of the proof. 

In the second part of the paper (\S7 - \S11), we study partitions with 
prescribed bounds on successive ranks and convert these into 
Rogers-Ramanujan partitions with weights. This leads to two new sets 
of partition theorems - three results involving partitions into parts 
$\not\equiv 0, \pm i$ ($mod$ 6), and three more concerning partitions 
into parts  $\not\equiv 0, \pm i$ ($mod$ 7), for $i=1,2,3$. In \S7 we 
describe the necessary preliminaries and historical background, namely, 
the work of Andrews and others [8], [10]-[14]. 
The proof of Theorem 2 connecting weighted Rogers-Ramanujan partitions 
with partitions into parts $\not\equiv 0, \pm 1$ ($mod$ 6) is given in full 
in \S8. Theorem 3 dealing with partitions into parts 
$\not\equiv 0, \pm 2$ ($mod$ 6) and Theorem 4 concerning partitions into 
parts $\not\equiv 0, \pm 3$ ($mod$ 6) are stated in \S9, but their proofs 
are only sketched because they are similar to the proof of Theorem 2. 
We draw attention here that the condition $\not\equiv 0, \pm 3$ ($mod$ 6) 
has to be interpreted properly, and when done so, leads to a very 
interesting new result (Theorem 4 of \S9) involving unrestricted partitions. 

In \S10 - \S11 we state  three new theorems connecting weighted 
Rogers-Ramanujan partitions with partitions into parts 
$\not\equiv 0, \pm i$ ($mod$ 7), $i=1,2,3$. Once again, only Theorem 5 which 
deals with parts $\not\equiv 0, \pm 3$ ($mod$ 7) is proved in full in \S10 
whereas proofs of Theorems 6 and 7 dealing with parts 
$\not\equiv 0, \pm i$ ($mod$ 7), $i=2,1$ are only sketched in \S11
because the details are similar. 

The theory of partitions is rich in examples of identities whose 
combinatorial interpretation yields the equality of partition functions 
defined in very different ways. Recently Alladi [2], [3], has developed 
a theory of weighted partition identities which deals with partition 
functions which are unequal, but where equality is attained by attaching 
weights. In [2], [3], the results primarily deal with the case where 
one set of partitions is a subset of the other, and where positive integral 
weights are attached to the smaller set of partitions.The results in this 
paper are further significant examples of such weighted identities and also 
deal with the more general situation of two unequal sets of partitions 
where weights positive or negative could be attached to either set. 
Theorem 1 of \S1 is of great interest because it is a reformulation of a 
deep theorem of G\"ollnitz [15] and has several important consequences 
(see \S2 - \S6). Theorems 2 - 7 (see \S7 - \S11) are appealing because the 
Rogers-Ramanujan partitions which had traditionally been associated with 
the modulus 5 (see Andrews [9; Ch. 7]), are for the first time, by means of 
weights, connected to 
partitions which are defined using congruence conditions modulo 6 and 7. 

Finally, in \S12, we describe briefly some future prospects of this research.

We conclude this section by describing some notation.

We denote the set of all Rogers-Ramanujan partitions by $\r$, and the 
subset of $\r$ consisting of partitions not having 1 as a part by $\r_2$. 
Given $\pi\in \r$, by a chain we mean a maximal block of 
parts differing by 2. Thus every $\pi \in \r$ can be 
decomposed into chains, and the parts in a given chain are all of the 
same parity.  For a chain $\chi$, we define its length $\ell=\ell(\chi)$ to
be the number of parts in $\chi$, and $\lambda = \lambda (\chi)$ to be its
least part.  

For any partition $\pi$, by $\sigma(\pi)$ we mean the sum of the 
parts of $\pi$, and by $\nu(\pi)$, the number of parts of $\pi$. We also 
use the standard notation that for any complex number $a$, the symbol 
$(a)_n$ is defined by
$$
(a)_n=(a;q)_n=\prod_{j=0}^{n-1}(1-aq^j)
$$
and 
$$
(a)_\infty=(a;q)_\infty={\lim_{n\to\infty}}(a)_n, \qquad \text{for}\quad |q|<1.
$$
As can be seen from the above, when the base in the product is $q$, it is 
often suppressed, but not when it is anything other than $q$. Further notation 
will be introduced when necessary. 

\newpage

\head \S1. {\un{A reformulation of the G\"ollnitz theorem}}\endhead

Given a Rogers-Ramanujan partition $\pi$, decompose it into chains. 
For a chain $\chi$ of length $\ell$, we define its weight $\omega(\chi)$ by 
$$
\omega(\chi)= \cases c^{\ell -1}(c+a b), \text{ if } \lambda (\chi) \text{ is
even },\\ 
a^{\ell} + (1+c)\sum_{k=1}^{{\ell}-1} a^k b^{{\ell-k}} + b^{\ell}, 
\text{ if } \lambda (\chi)=1,\\  
(1+c){\left\{a^{\ell} +(1+c)\sum_{k=1}^{{\ell}-1} a^k b^{{\ell-k}} + b^{\ell}
+ b^{\ell}\right\}}, \text{ if } \lambda (\chi) >1 \text{ is odd },\endcases
\tag1.1
$$
where $a, b, \text{ and }c$ are free parameters whose role
will be described soon.  Finally the weight $\omega(\pi)$ of a
Rogers-Ramanujan partition is defined as the product of the weight of
its chains $\chi$; that is
$$
\omega(\pi) = \prod_{\chi}\omega(\chi).\tag1.2
$$

Next, let ${\overset\rightarrow\to\pi}= (\pi_1; \pi_2; \pi_3)$
denote a vector partition with $\pi_1$ and $\pi_2$ having distinct odd
parts, and $\pi_3$ having distinct even parts.  Denote by $\V$ the
set of all such vector partitions.

Our goal is prove the following result and discuss its implications.

\proclaim{Theorem 1} Let $\r, \V \text{ and }\omega$, be as above.  Then
for any integer $n\ge 0,$ 
$$\sum_{\pi\in\r, \sigma(\pi)=n}\omega(\pi)=
\sum_{{\overset\rightarrow\to\pi}\in\V, 
\sigma(\overset\rightarrow\to\pi)=n}a^{\nu(\pi_1)} b^{\nu(\pi_2)}
c^{\nu(\pi_3)}.
$$
\endproclaim

The proof of Theorem 1 is given at the end in $\S 6.$  It is based on
the method of weighted words due to Alladi-Andrews-Gordon [5] which
provided a refinement and generalization of a deep theorem of
G\"ollnitz [15].  The main ideas of [5] are described in $\S 5$ as the
necessary background for the proof of Theorem 1 which is given in $\S6.$

Before we give the proof of Theorem 1, we discuss its implications in 
the next three sections.

\head \S2. {\un{Jacobi's triple product identity}}\endhead

From the definition of $\V$ it follows that
$$
(-aq;q^2)_\infty(-bq;q^2)_\infty(-cq^2;q^2)_\infty
=\sum_{n\ge 0}\sum\Sb
\overset\rightarrow\to\pi \in\V \\ 
\sigma(\overset\rightarrow\to\pi)=n\endSb 
a^{\nu(\pi_1)}b^{\nu(\pi_2)} c^{\nu(\pi_3)}q^n.\tag2.1
$$
Now take
$$
ab=1, \quad\text{and}\quad c= -1.\tag2.2
$$
So the product on the left in (2.1) is
$$
(-aq;q^2)_\infty(-a^{-1}q;q^2)_\infty(q^2;q^2)_\infty.\tag2.3
$$
The choices (2.2) imply that
$$
c+ab =0 \quad\text{ and }\quad 1+c=0.\tag2.4
$$
Therefore for partitions $\pi \in\r,$ the only chains with non-zero
weights will be the chains
$$\chi: 1+3+5+...+(2n-1), n\ge1,\tag2.5$$
for which
$$\omega(\chi)=a^n+b^n=a^n+a^{-n}, n\ge 1.\tag2.6$$
Thus the only partitions $\pi \in \r$ which will have non-zero weights
will be the partitions of $n^2$ given by $1+3+...+(2n-1)$ with weights
$\omega(\pi)=\omega(\chi)$ as in (2.6). Thus by Theorem 1 and (2.2) 
through (2.6) we get
$$
1+\sum_{n=1}^{\infty}(a^n+a^{-n})q^{n^2}=\sum_{n=-\infty}^{\infty}a^nq^{n^2}
=(-aq;q^2)_\infty(-a^{-1}q;q^2)_\infty(q^2;q^2)_\infty,\tag2.7
$$
which is Jacobi's triple product identity for theta functions.

\head \S3. {\un{Sylvester's refinement of Euler's theorem}}\endhead

In an important paper of 1882, Sylvester [17] improved many
partition results of Euler by exploiting combinatorial techniques.  In
particular, Sylvester proved:

\proclaim{Theorem S:}  Let $k$ and $n$ be positive integers.  Then the
number of partitions of $n$ into odd parts of which exactly $k$ are
different, equals the number of partitions of $n$ into distinct parts
which can be grouped into $k$ maximal blocks of consecutive integers.
\endproclaim

Euler's famous theorem on the equality of partitions of $n$ into odd
parts and distinct parts follows from Theorem S by summing over $k$.

We now show that Theorem S follows from Theorem 1.  To this end, take
$$
c=1, \quad a+b=0.\tag3.1
$$
Then the product in (2.1) can be rewritten as
$$
\overset\infty\to{\underset m=1\to\prod}(1+abq^{4m-2})(1+q^{2m})=
\overset\infty\to{\underset m=1\to\prod}\left
(\dfrac{1+abq^{4m-2}}{1-q^{4m-2}}\right ).\tag3.2
$$
Note that since
$$
\dfrac{1+xq^n}{1-q^n}= 1+(1+x)\left(q^n+q^{2n}+q^{3n}+...\right),\tag3.3
$$
the product on the right in (3.2) has the interpretation that it is
the generating function for partitions into parts $\equiv$ 2(mod 4)
where such partitions $\pi$ are counted with weights
$$
(1+ab)^{\nu_d (\pi)},\tag3.4
$$
where $\nu_d(\pi)$ is the number of different parts of $\pi.$

Next observe that with $c=1$ we have
$$
a^n+2\sum_{j=1}^{n-1}a^{n-j}b^j+b^n=(a+b)\sum_{j=0}^{n-1}a^{n-1-j}b^j.\tag3.5
$$
So if $(a+b)=0$, then it follows from (3.5) that all odd chains in
(1.1) must have weight 0.  Thus the only chains with non-zero weights
are the partitions into even parts differing by $\ge 2.$ The weight of
such partitions given by (1.2) will be 
$$(1+ab)^k,\tag3.6$$
where $k$ is the number of even chains. Thus from (3.4) and
(3.6), Theorem S follows dilated by a factor of 2, where the odd
numbers are replaced by integers $\equiv 2 (\mod  4)$ and the positive
integers are replaced by even numbers.

\head \S4. {\un{Two weighted partition theorems}}\endhead

Take
$$
b=c=1, \text{ and } a=0,\tag4.1
$$
in Theorem 1.  Then the product (2.1) is
$$
\overset\infty\to{\underset m=1\to\prod}(1+q^m)\quad =
\quad\sum_{n=0}^{\infty} D(n) q^n,\tag4.2
$$
the generating function of $D(n),$ the number of partitions of $n$
into distinct parts.

Next, the choices in (4.1) imply that in (1.1) all even chains have
weight 1, all odd chains $\chi$ with $\lambda(\chi)=1$ have weight 1, 
and all odd chains $\chi$ with $\lambda(\chi)> 1$ have weight 2.  
Thus by (1.2) 
$$
\omega(\pi)=2^k,\tag4.3
$$
where $k$ is the number of odd chains $\chi$ of $\pi$ with $\lambda(\chi)>1.$
So we get the following result as a consequence of Theorem 1:

\proclaim{Theorem A} (Alladi [2])Let $D(n)$ denote the number of partitions of
$n$ into distinct parts.  Given $\pi\in \r,$ let its weight
$\omega(\pi)=2^k,$ where $k$ is the number of odd chains $\chi$ of $\pi$ with
$\lambda(\chi)>1.$  Then
$$
\sum\Sb\pi\in\r \\ \sigma(\pi)=n \endSb \omega(\pi) = D (n).
$$
\endproclaim

\un{Remarks}: Theorem 10 of [2] is actually the same result as Theorem A, 
but is stated differently.  In [2], $k$ is interpreted as the number of gaps
between the odd parts of $\pi$ and $-1$ which are $>2.$  The proof of
this result in [2] involves constructing a (combinatorial) surjective
map between the set of partitions into distinct parts and its subset,
namely the set of Rogers-Ramanujan partitions.

Since partitions into odd parts are equinumerous with partitions into
distinct parts, it is interesting to ask whether there is another
weighted partition identity connecting partitions into odd parts and
Rogers-Ramanujan partitions.  For this purpose take
$$
a=b=c=1\tag4.4
$$
in Theorem 1.  Now rewrite the product in (2.1) as
$$
(-q;q^2)_\infty(-q;q^2)_\infty(-q^2;q^2)_\infty= 
\overset\infty\to{\underset m=1\to\prod}(1+q^{2m-1})(1+q^m)=
\overset\infty\to{\underset m=1\to\prod}\left(\dfrac{1+q^{2m-1}}{1-q^{2m-1}}
\right).\tag4.5
$$
By (3.3), the product in (4.5) has the interpretation that it is the
generating function of partitions in $\pi'$ into odd parts, where
$\pi'$ is counted with weight
$$
\omega_1(\pi')=2^{\nu_d(\pi')}.\tag4.6$$
In (4.6), as in Sylvester's theorem, $\nu_d(\pi')$ is the number of
different parts of $\pi'$.

Finally, note that with the choices in (4.4), the weights in (1.1) become
$$
\omega(\chi)=\cases 2, \text{ if } \chi\text{ has even parts } \\ 2\ell,
\text{ if } \ell(\chi)=\ell \text{ and }\lambda (\chi)=1,\\ 4\ell, \text{ if
} \ell (\chi)=\ell, \text{ and } \lambda (\chi) > 1.\endcases\tag4.7
$$
With these values of $w(\chi),$ and with $w(\pi)$ defined by (1.2), we
get from (4.4), (4.5), and (4.6) the following new weighted partition
theorem:
\proclaim{Theorem B} Let $\Cal O$ denote the set of partitions into
odd parts.  Then
$$
\sum_{\pi\in\r, \sigma(\pi)=n}\omega(\pi)=
\sum_{\pi'\in\Cal O, \sigma(\pi')=n}2^{\nu_d(\pi')}.
$$
\endproclaim

{\un{Remarks}}: In \S8, by replacing $D(n)$ with the number of partitions of 
$n$ into odd parts, and by thinking of odd parts as being in residue classes 
$\equiv1,3,5(mod$ $6)$, Theorem A is reformulated and a new proof given 
(see Theorem 3). This has the advantage of producing two other similar 
results (Theorem 2 of \S8 and Theorem 4 of \S9), both of which are new.

\head \S5. {\un{the method of weighted words}}\endhead

In 1967,  G\"ollnitz [15] proved the following deep partition theorem:

\proclaim{Theorem G} Let $B(n)$ denote the number of partitions of
$n$ into distinct parts $\equiv 2, 4, \text { or } 5(\mod 6).$

Let $C(n)$ denote the number of partitions of $n$ in the form
$m_1+m_2\cdots + m_{\nu}$ such that $m_{\nu}\ne 1,3$, and $m_i-m_{i+1}\ge 6$ 
with strict inequality if $m_i\equiv$ 0, 1, or 3(mod 6). Then
$$
B(n)=C(n).
$$
\endproclaim

In [5], Alladi, Andrews, and Gordon, obtained substantial refinements
and generalizations of Theorem G by using a technique called the
method of weighted words.  We now describe briefly the main ideas in
[5].

Theorem G is viewed in [5] as emerging out of the {\it{key identity}}
$$
\sum_{i,j,k} a^i b^j c^k \sum\Sb i=\alpha+\delta +\varepsilon\\
j=\beta+\delta +\phi\\ k=\gamma+\varepsilon+\phi\endSb
\dfrac{q^{T_s+T_\delta+T_\varepsilon+T_{\phi-1}} (1
-q^\alpha(1-q^\phi))}{(q)_\alpha(q)_\beta(q)_\gamma(q)_\delta(q)_\varepsilon
(q)_\phi}=(-aq)_\infty(-bq)_\infty(-cq)_\infty\tag5.1 
$$
under the {\it{standard transformations}}
$$
\matrix\text{(dilation)}\quad q\mapsto q^6,\\
\text{(translations)}\quad a\mapsto aq^{-4}, b\mapsto bq^{-2}, c\mapsto
cq^{-1}.\endmatrix\bigg\}\tag5.2
$$
In (5.1), $s=\alpha+\beta+\gamma+\delta+\varepsilon+\phi$, $T_m=m(m+1)/2$. 
Clearly when the transformations (5.2) are applied to the product in
(5.1) we get 
$$
\prod_{m=1}^{\infty}(1+aq^{6m-4})(1+bq^{6m-2})(1+cq^{6m-1})\tag5.3
$$
which is the three parameter refined generating function of $B(n)$ in
Theorem G.  We now describe how the series in (5.1) becomes the
refined generating function of $C(n)$ under the influence of (5.2).

We consider the integer 1 as occurring in three primary colors $a$, $b$,
and $c$, and integers $n\ge 2$ as occurring in the three primary colors
as well as in three secondary colors $d=ab, e=ac$, and $f=bc.$
The integer $n$ in color $a$ is denoted by the symbol $a_n,$ with
similar interpretation for the symbols $b_n,...,f_n.$  In order to
discuss partitions (words) involving the symbols, we need an ordering
among them, and the one we choose is 
$$
a_1<b_1<c_1<d_2<e_2<a_2<f_2<b_2<c_2<d_3<e_3<a_3<f_3<\dots\tag5.4
$$
The reason for choosing this ordering is because under the transformations 
(5.2) the symbols become
$$
\matrix a_n\mapsto 6n-4, b_n\mapsto 6n-2, c_n\mapsto6n-1,\text { for } 
n\ge 1, \\
d_n=ab_n\mapsto 6n-6, e_n=ac_n\mapsto 6n-5, f_n=bc_n\mapsto 6n-3,
\text{ for }n\ge 2,\endmatrix\bigg\}\tag5.5
$$ 
and so (5.4) becomes
$$
2<4<5<6<7<8<9<10<11<12<13<\dots,\tag5.6
$$
the natural ordering among the integers.

It is convenient to write down the complete list of symbols in (5.4), namely,
$$
\und{e}_1<a_1<\und{f}_1 <b_1 <c_1<d_2<e_2<a_2<f_2<b_2<c_2<d_3<\dots\tag5.7
$$
where $\und{e}_1$ and $\und{f}_1$ in (5.7) are underlined because they
do not really occur.  We have omitted writing $\und{d}_1$ in
(5.7) because $\und{d}_1=0$.

Next, let $x_n$ denote the symbol occupying position $n$ in (5.7); 
that is $x_1=\und{e}_1, x_2=a_1, x_3=\und{f}_1,$ 
and so on.  By a Type-1 partition we mean an expression of the form
$x_{n_1}+ x_{n_2}+ \dots + x_{n_\nu}$, where the
$x_{n_i}$ are chosen from the non-underlined set in (5.7) and
satisfy the {\it{standard gap conditions}}
$$
n_i-n_{i+1}\ge 6,\text{ with strict inequality if } x_{n_i} 
\text{ is of secondary color.}\tag5.8
$$
The main result in [5] is
\proclaim{\un{Theorem C}} Let $\vee(n; i, j, k)$ denote the number of
vector partitions $(\pi'_1, \pi'_2, \pi'_3)$ of $n$ such that $\pi'_1$
has $i$ distinct parts all in color $a$, $\pi'_2$ had $j$ distinct parts 
all in color $b$, and $\pi'_3$ has $k$ distinct parts all in color $c$.

Let $C(n;\alpha, \beta, \gamma, \delta, \varepsilon, \phi)$ denote the
number of Type-1 partitions of $n$ having $\alpha$ a-parts, $\beta$
b-parts, $\dots, \phi$ f-parts.  Then
$$
\vee(n;i,j,k)=\sum\Sb i=\alpha +\delta+\varepsilon\\
j=\beta+\delta+\phi\\k=\gamma +\varepsilon+\phi\endSb C(n;\alpha,
\beta, \gamma, \delta, \varepsilon,\phi).
$$
\endproclaim

Clearly the generating function of $\vee(n;i,j,k)$ is
$$
\sum_{i,j,k,n} \vee(n;i,j,k)a^ib^ic^kq^n=
(-aq)_\infty(-bq)_\infty(-cq)_\infty.\tag5.9
$$
In [5] it is shown that for given $\alpha,\beta,\gamma,\delta,\varepsilon, 
\text{ and }\phi,$
$$
\sum_n C(n;\alpha, \beta, \gamma, \delta, \varepsilon, \phi)q^n=
\dfrac{q^{T_s+T_\delta+T_\varepsilon+T_{\phi-1}}(1-q^\alpha(1-q^\phi))}
{(q)_\alpha(q)_\beta(q)_\gamma(q)_\delta(q)_\varepsilon(q)_\phi}.\tag5.10
$$
Thus Theorem 3 is a consequence of (5.9), (5.10) and the key identity (5.1).  

In this approach, under the transformations (5.1), the primary colors
$a, b, c,$ correspond to the residues 2, 4, 5 (mod 6) which determine
the partition function $B(n)$ in Theorem G.  Thus the secondary colors
are 2+4$\equiv$6(mod 6), 2+5$\equiv$7(mod 6), and 4+5$\equiv$9(mod 6), 
meaning that the residue classes 0, 1, 3 (mod 6),
represent secondary colors, but parts in these residue classes are all
$\ge 6$ because $\und{d}_1, \und{e}_1, \und{f}_1,$ do not occur.  This
explains the condition $m_{\nu} \ne 1 \text{ or } 3$ in 
defining $C(n)$ in Theorem G.  Also the strict inequality $m_i-m_{i+1}
> 6$ when $m_i\equiv$ 0,1, or 3(mod 6) is to be interpreted as the
inequality being strict when $m_i$ is of secondary color.  The
difference conditions defining Type-1 partitions translate to those
defining $C(n)$ in Theorem G when the standard transformations (5.2)
are applied.  Thus Theorem C is a strong refinement and generalization
of Theorem G.

\head \S6. {\un{Proof of Theorem 1}}\endhead

Instead of the standard transformations (5.2), let us consider the
effect of the {\it{quadratic transformations}}
$$
\matrix\text{(dilation)}\quad q\mapsto q^2,\\
\text{(translations)}\quad a\mapsto aq^{-1}, b\mapsto bq^{-1}, c\mapsto c,
\endmatrix\bigg\}\tag6.1
$$
on (5.1). Clearly (6.1) converts the product on the right in (5.1) to
the product in (2.1) which is the generating function of partitions
$\overset\rightarrow\to\pi \in\V$ in Theorem 1.

The effect of (6.1) on the symbols is
$$
\matrix a_n\mapsto 2n-1,\quad b_n\mapsto 2n-1, \quad c_n\mapsto 2n,
\quad\text{ for }n\ge1,\\d_n=ab_n\mapsto 2n-2, e_n=ac_n\mapsto 2n-1,
f_n=bc_n\mapsto 2n-1, \text{ for } n\ge 2.
\endmatrix\bigg\}\tag6.2
$$
Thus (5.4) becomes
$$
1_a<1_b<2_c<2_{ab}<3_{ac}<3_a<3_{bc}<3_b<4_c<4_{ab}<5_{ac}<5_a<\dots,\tag6.3
$$
where we have reversed the convention by indicating the color with a
subscript.  In this case Type-1 partitions are Rogers-Ramanujan
partitions satisfying certain color conditions which determine the 
weight of such a partition.  These conditions
imply that when two integers in (6.3) differ by $>2$, then colors can
be attached to either integer in all possible ways.  That is there is
no interference or dependence here.  What this means is that in order
to determine the weight of a Rogers-Ramanujan partition, we need only consider
chains of parts and calculate the weights of these chains; then
by the independence, the weight of the partition can be calculated
using the product formula (1.2).

Consider now a chain $\chi$ of even parts $m_1>m_2\dots >m_\ell.$  Note
that (5.8) implies that we have a choice only for $m_1$ to have
color $c\text{ or } ab,$ but all parts $,m_1$ in the chain must
have color $c$.  Thus in this case
$$\omega(\chi)
=c^{\ell-1} (c+ab), \text{ if } \lambda(\chi) \text{ is even},\tag6.4
$$
as in (1.1).

Next consider a chain $\chi$ with $\lambda(\chi)=1.$  If any part in this 
chain has color $b$, then the next higher part (and therefore all parts
higher) must have color $b$.  If any part in this chain has color
$bc$, then all parts higher must have color
$b$.  If any part in this chain has color $a$, then the next higher
part can have color $a, bc, \text{ or } b.$  Note that no part in this
chain can have color $ac$ since the chain has to start with colors $a
\text{ or } b.$  So we have the following cases to consider.

\un{Case 1}: $1 \text{ has color }b$.

Then all parts have color $b$. So the weight of this chain is $b^n$.

\un{Case 2}:  All parts have color $a.$

Then clearly the weight of the chain is $a^n.$

\un{Case 3}:  The only parts in color $a$ are $1, 3, \dots, 3k-1,$
with $1\le k<\ell=\ell(\chi).$  

Then $2k+1$ can have color $bc$ or $b$, but all parts $>2k+1$ must have color 
$b$.  So the weight of this chain is
$$
a^k(bc+b) b^{\ell-k-1}=(1+c)a^k b^{\ell-k}.
$$
We need to sum this over all $k$ to get the weight of all chains covered by
Case 3.

The sum of the weights of chains in Cases 1, 2, and 3 is
$$
b^\ell + a^\ell + (1+c) \sum_{k=1}^{{\ell}-1} a^k b^{{\ell-k}},\tag6.5
$$
as in (1.1).

Finally, consider chains $\chi$ with $\lambda (\chi)$=odd$>1$. In
this situation we have color choices as in Cases 1, 2, 3 above, plus
the cases where $\lambda(\chi)$ has color $ac\text{ or } bc.$  So this gives
rise to three more cases.

\un{Case 4}: $\lambda(\chi)$ has color $ac,$ and next part has color $bc.$

Then the rest of the parts have color $b.$  So the weight of this chain is 
$$
(ac)(bc) b^{\ell-2} = c^2 ab^{\ell-1}.\tag6.6
$$

\un{Case 5}: $\lambda (\chi)$ has color $ac$, and next part has color $a$
or $b$.  

So the chain after $\lambda(\chi)$ has length $\ell-1$ and this
situation covered by the colorings as in Cases 1, 2, and 3.  So the
weight in this case is 
$$
\matrix ac\left\{a^{\ell-1} +(1+c)(a^{\ell-2} b+a^{\ell-3} b^2+
\dots + ab^{\ell-2})+b^{\ell-1}\right\}\\ = a^\ell c+
c(1+c)(a^{\ell-1}b+a^{\ell-2}b^2+\dots + a^2 b^{\ell-2})+c
ab^{\ell-1}.\endmatrix\tag6.7
$$

\un{Case 6}:\quad $\lambda(\chi)$ has color $bc$.

Then the rest of the parts all have color $b$.  So the weight of this
chain is 
$$
(bc) b^{\ell-1} =cb^\ell.\tag6.8
$$

So adding the weights in (6.6), (6.7) and (6.8) we get 
$$c a^\ell + c(1+c)(a^{\ell-1} b+a^{\ell-2} b^2+\dots + ab^{\ell-1}) +
cb^\ell\tag6.9$$
Finally we need to add the weights in (6.9) and (6.5) to get the
weights of chains in Cases 1 through 6 that cover all chains with
$\lambda(\chi)=\text{ odd }>1.$  This gives
$$
(1+c)a^\ell+(1+c)^2(a^{\ell-1} b+a^{\ell-2} b^2 + \dots + ab^{\ell-1})
+(1+c) a^\ell\tag6.10
$$
which is what is given in (1.1).

Thus the weights in (1.1) have been established and this completes the
proof of Theorem 1.

{\un{Remarks}}: Previously we had discussed consequences of (5.1) and 
Theorem 1 under the dilations $q\mapsto q^3$ (see [1]) and $q\mapsto q^4$ 
(see [4]), and certain sets of translations. These cubic and quartic 
transformations lead to different combinatorial versions of Theorem G.  

\head \S7. {\un{Successive ranks with prescribed bounds}}\endhead

The Ferrers graph of every partition contains a Durfee square, namely, 
the largest square of nodes starting from the top left hand corner of 
the graph. Through every node on the descending diagonal of the 
Durfee square there is a hook passing through it, namely, the set of nodes 
from (and including) that node on the diagonal horizontally to its right 
and vertically below it. The Ferrers graph of the partition 7+6+6+4+4+2+1+1, 
its Durfee square, and its hooks, are illustrated below.

\line{\hfill\epsfxsize=5truein\epsfbox{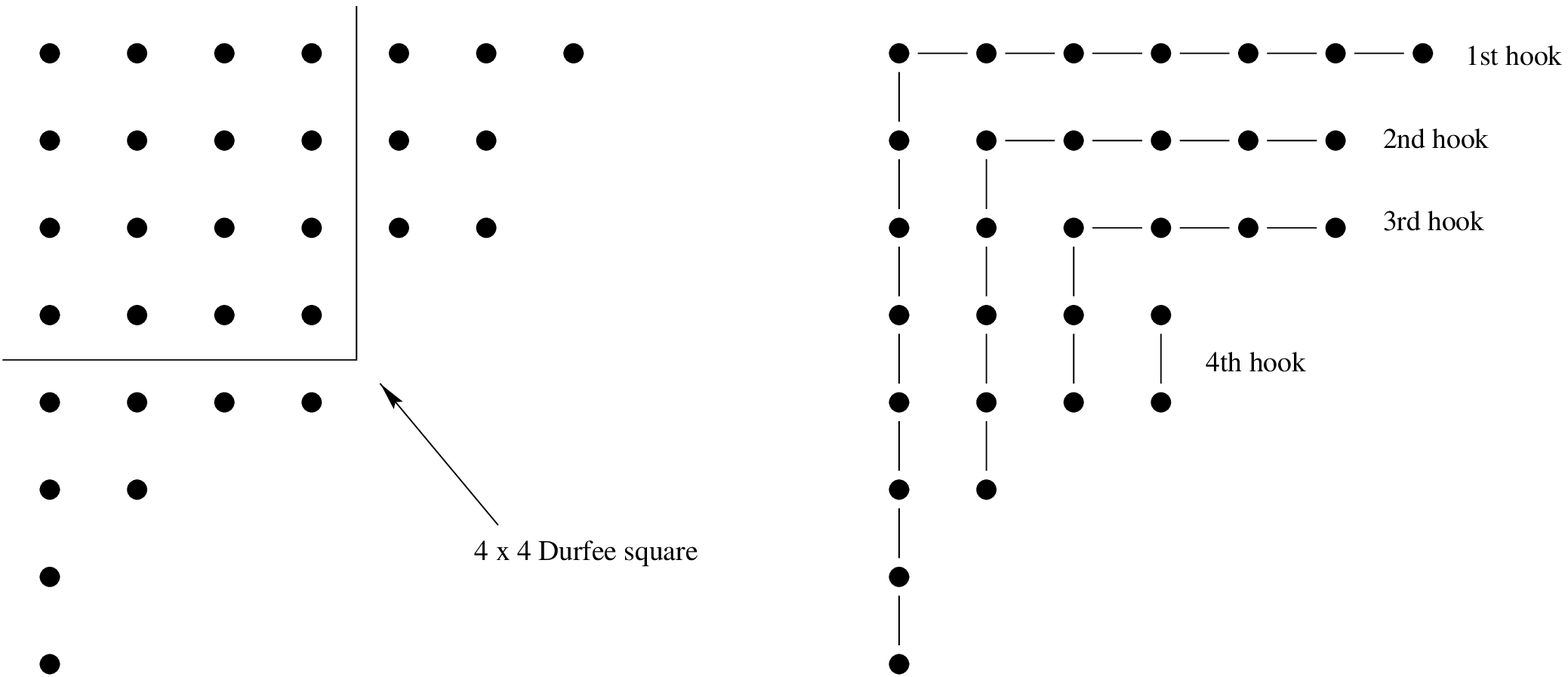}\hfill}
\centerline{\smallcaps Figure I}

If the graph of a partition has a $k\times k$ Durfee square, then the 
number of hooks is k. The largest hook is called the first hook, 
the second largest hook is called the second hook, and so on. 

Given a hook of a Ferrers graph, its rank is the number of horizontal nodes 
to the right of the vertex of the hook, minus the number of vertical nodes 
below the vertex of the hook. The rank of the i-th hook is called the i-th 
rank of the partition. The successive ranks in Figure 1 are -1, 0, 1, and -1. 

The first rank is, of course, the familiar rank of the partition made famous 
by Dyson [14] in conjecturing a combinatorial explanation of Ramanujan's 
congruences modulo 5 and 7 for the partition function. The Dyson conjectures 
involving the rank were proved by Atkin and Swinnerton Dyer [12]. Atkin [11] 
was led to consider successive ranks while attempting a study of 
Ramanujan's partition congruence mod 11. Following Atkin, Andrews [8], and 
Bressoud [13] studied partitions with prescribed bounds for successive 
ranks. In [10] the notion of successive ranks was generalized 
to hook differences, with the hook vertices not necessarily on the diagonal 
of the Durfee square. Also in [10] the succesive rank theorem of Andrews [8] 
and Bressoud [13] was revised as follows:

\newpage

\proclaim{Theorem R} Let $Q_{k,i}(n)$ denote the number of partitions of 
$n$ such that the successive ranks all take values in the interval 
$[-i+2, k-i-2]$. 

Let $A_{k,i}(n)$ denote the number of partitions into parts 
$\not\equiv 0, \pm i$ ($mod$ $k$). 

Then for $1\le i< k/2$ we have
$$
Q_{k,i}(n)=A_{k,i}(n).
$$
\endproclaim
  
Although the case $2i=k$ is not covered by Theorem R, it is possible to 
deal with this case by defining $A_{k,i}(n)$ via the identity
$$
\sum_{n=0}^{\infty}A_{k,i}(n)q^n=
\frac{(q^k;q^k)_\infty(q^i;q^k)_\infty(q^{k-i};q^k)_\infty}{(q)_\infty}.
\tag7.1
$$

If $2i\ne k$, then $A_{k,i}(n)$ defined in (7.1) has the partition 
interpretation as in Theorem R. However, if $2i=k$, then $A_{k,i}(n)$ 
does not have a partition interpretation in the standard sense. The advantage 
of (7.1) is that it leads to 

\proclaim{Theorem R'}
$$    
Q_{2k,k}(n)=A_{2k,k}(n).
$$
\endproclaim
\bigskip

\un{Remarks}: Theorem R' is a consequence of Theorem 5 of [10], but we stress 
here that if $A_{2k,k}(n)$ is to be interpreted as the number of 
partitions of $n$ into parts $\not\equiv 0, \pm k$ ($mod$ $2k$), then 
this has to be in the sense of (7.1), where the 
residue class $k(mod$ $2k)$ is ``deleted twice'' because it occurs as 
both $k$ and $-k$ $(mod$ $2k).$ 

\head \S8. {\un{A weighted partition theorem mod 6}}\endhead
 
Given a partition $\pi$ whose Ferrers graph has a $k\times k$ Durfee square, 
consider the partition $\rho (\pi)=\overset\sim \to \pi$ into $k$ parts 
$h_1+h_2+...+h_k$, where $h_i$ is the number of nodes in (= the length of) 
the i-th hook of $\pi$. Note that $h_i-h_{i+1}\ge 2$ for $1\le i\le k-1$ 
and so $\overset\sim \to \pi$ is a Rogers-Ramanujan partition. The mapping 
$$
\pi\rightarrow\rho (\pi)=\overset\sim \to \pi\tag8.1
$$
is surjective. Also, if $r_i(\pi)$ is the i-th rank of $\pi$, then under this 
mapping
$$
h_i-r_i(\pi)\equiv1(mod \quad 2).  \tag8.2
$$ 
The study of the surjection (8.1) along with Theorems R and R' 
will lead us to several weighted partition identities in this and subsequent 
sections. We begin with

\proclaim{Theorem 2} Given $\overset\sim \to \pi \in \r_2$, let its weight 
$\omega_2(\overset\sim \to \pi)=2^r$, where $r$ is the number of even chains 
$\chi$ in $\overset\sim \to \pi$ with $\lambda(\chi)>2$. 

Let $A_{6,1}(n)$ denote the number of partitions of $n$ into parts 
$\not\equiv 0, \pm 1$ ($mod$ $6$). Then 
$$
A_{6,1}(n)=\sum_{\overset\sim \to \pi\in\r_2, \sigma(\overset\sim \to \pi)=n}
\omega_2(\overset\sim \to \pi).
$$
\endproclaim

{\bf{Proof}}: Take $k=6, i=1$ in Theorem R. Thus 
$$
A_{6,1}(n)=Q_{6,1}(n).\tag8.3
$$
We will show that
$$
Q_{6,1}(n)=\sum_{\overset\sim \to \pi\in \r_2, \sigma(\overset\sim \to \pi)=n}
\omega_2(\overset\sim \to \pi).\tag8.4
$$
Theorem 2 will follow from (8.3) and (8.4). 

We know that $Q_{6,1}(n)$ is the number of partitions of $n$ whose 
successive ranks take values 1,2, or 3. Given a partition $\pi$ enumerated 
by $Q_{6,1}(n)$, consider the partition $\rho(\pi)=\overset\sim \to \pi$ 
generated by the hooks of $\pi$. If the i-th part of $\overset\sim \to \pi$ 
(= i-th hook length of $\pi$) is odd, then by (8.2) the i-th rank of 
$\pi$ must be even, and so must be 2. So there is only one way in which this 
hook can occur. For instance, if 7 is the hook length and the rank is 2, 
then the hook must be

\bs
\line{\hfill\epsfxsize=1truein\epsfbox{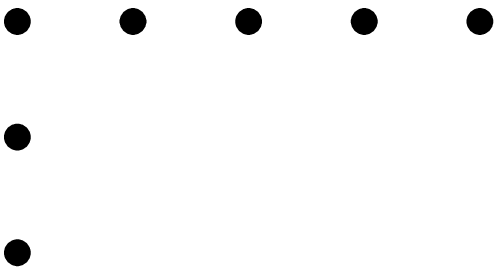}\hfill}
\centerline{\smallcaps Figure II}
\bs

Note that 1 can never occur as a hook length with rank 2. Thus 1 cannot 
occur as a part of $\overset\sim \to \pi$. Hence 
$\overset\sim \to \pi \in \r_2$. 

Now if the i-th part of $\overset\sim \to \pi$ is even, then the i-th rank of 
$\pi$ could be either 1 or 3. The question is under what circumstances can 
both values 1 and 3 occur as the rank? 

The integer 2 as a hook length can occur only with rank 1 and cannot have 
rank 3. But even integers $>2$ can arise as hook 
lengths in two ways, one with rank 1, and another with rank 3. For instance, 
6 as a hook length can be realized as 

\bs
\line{\hfill\epsfxsize=3truein\epsfbox{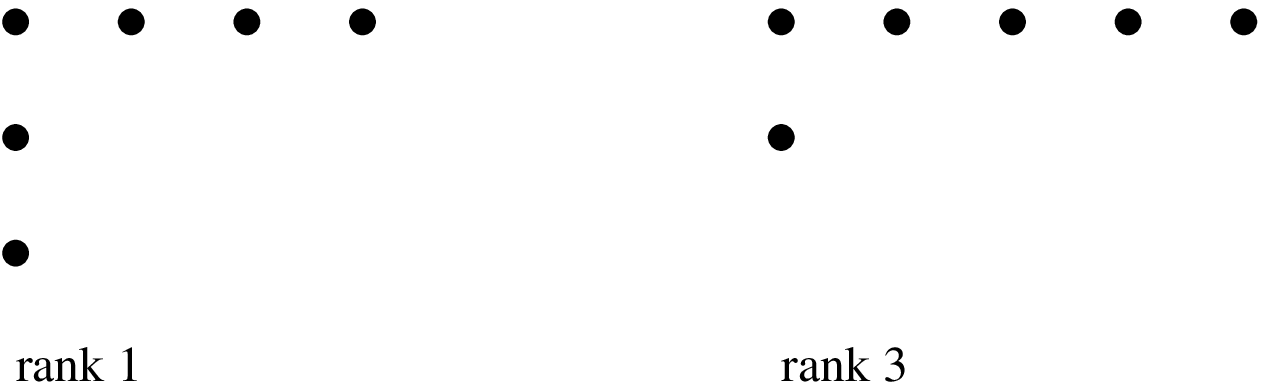}\hfill}
\centerline{\smallcaps Figure III}
\bs

  If two consecutive integers $2j$ and $2j+2$ occur as parts of 
$\overset\sim \to \pi$, and if the hook representing $2j$ has a certain rank, 
then the hook representing $2j+2$ must have the same rank. On the other 
hand, if two consecutive parts of $\pi$ differ by $>2$, then the rank 
of one part is independent of the rank of the other. Because of this 
independence, the weight to be attached to $\overset\sim \to \pi$ can 
be computed by decomposing $\overset\sim \to \pi$ into chains and taking 
the product of the weights of chains.  

So given $\overset\sim \to \pi \in \r _2$, decompose it into chains. 
All odd parts of $\overset\sim \to \pi$ have unique representations as hook   
lengths of $\pi$. With regard to even chains $\chi$ of $\overset\sim \to \pi$, 
the hook representation of $\lambda(\chi)$ in $\pi$ uniquely determines 
the hook representation of all other parts in that chain. If $\lambda(\chi)
>2$ is even, then $\lambda(\chi)$ admits two representations as a hook, 
one with rank 1 and another with rank 3. If $\lambda(\chi)=2$, then the 
hook representation must have rank 1. Thus to each 
$\overset\sim \to \pi \in \r_2$, there corresponds $2^r$ partitions 
partitions $\pi$ enumerated by 
$Q_{6,1}$, such that $\rho(\pi)=\overset\sim \to \pi$. This establishes 
(8.4) and completes the proof of Theorem 2.

\newpage 
 
\head \S9. {\un{Two more weighted identities mod 6}}\endhead

Pursuing the ideas of \S8, we get two more results.

\proclaim{Theorem 3} Given a Rogers-Ramanujan partition 
$\overset\sim \to \pi$, let $\omega_3(\overset\sim \to \pi)=2^r$, where 
$r$ is the number of odd chains $\chi$ of $\overset\sim \to \pi$ with 
$\lambda(\chi)>1$. 

Let $A_{6,2}(n)$ denote the number of partitions of $n$ into parts 
$\not\equiv 0\pm2$ ($mod$ $6$). Then  
$$
A_{6,2}(n)=\sum_{\overset\sim \to \pi\in \r, \sigma(\overset\sim \to \pi)=n}
\omega_3(\overset\sim \to \pi).
$$
\endproclaim

{\bf{Proof}}: We know from Theorem R that 
$$
A_{6,2}(n)=Q_{6,2}(n).\tag9.1
$$
Recall that $Q_{6,2}(n)$ is the number of partitions $\pi$ of $n$ such that 
the successive ranks take values 0, 1, or 2.

Given $\overset\sim \to \pi \in \r$, consider a partition $\pi$ enumerated 
by $Q_{6,2}$ with $\rho(\pi)=\overset\sim \to \pi$. In this case, the 
even parts of $\overset\sim \to \pi$ have unique representations as 
hooks of $\pi$ because the rank must be 1. The odd parts of 
$\overset\sim \to \pi$ which are $>1$ can have two possible representations 
as hooks of $\pi$ - one with rank 0 and another with rank 2. In a chain 
$\chi$ of odd parts of $\overset\sim \to \pi$, the hook representation of 
$\lambda(\chi)$ uniquely represents the hook representation of all other 
parts in the chain. When $\lambda(\chi)>1$ is odd, it admits two hook 
representations in $\pi$, but $\lambda(\chi)=1$ admits only one representation 
with rank 0. Thus to each $\overset\sim \to \pi \in \r$, there 
corresponds $\omega_3(\overset\sim \to \pi)$ partitions $\pi$ enumerated 
by $Q_{6,2}$ such that $\rho(\pi)=\overset\sim \to \pi$. Therefore  
$$
Q_{6,2}(n)=\sum_{\overset\sim \to \pi\in \r, \sigma(\overset\sim \to \pi)=n}
\omega_3(\overset\sim \to \pi).\tag9.2
$$
Theorem 3 follows from (9.1) and (9.2). 

{\un{Remarks}}: Note that $A_{6,2}(n)$ is the number of partitions of $n$ 
into odd parts because 
$$
\sum_{n=0}^{\infty}Q_{6,2}(n)q^n=
\frac{1}{(q;q^6)_\infty(q^3;q^6)_\infty(q^5;q^6)_\infty}=
\frac{1}{(q;q^2)_\infty}.\tag9.3
$$
By Euler's theorem, $Q_{6,2}(n)$ equals $D(n)$, the number of partitions of 
$n$ into distinct parts. Thus Theorem 3 is a reformulation of Theorem A but 
here the result is proved using successive ranks. In doing so, we see a 
similarity with Theorem 2. 

We now consider the one remaining product mod 6, namely, 
$$
\sum_{n=0}^{\infty}A_{6,3}(n)q^n=
\frac{(q^6;q^6)_\infty(q^3;q^6)_\infty(q^3;q^6)_\infty}{(q)_\infty}=
$$
$$
\frac{(q^3;q^6)_\infty}{(q;q^3)_\infty(q^2;q^3)_\infty}=
\frac{1}{(q;q^3)_\infty(q^2;q^3)_\infty(-q^3;q^3)_\infty}.\tag9.4
$$
The product on the right in (9.4) is the generating function for 
unrestricted partitions $\pi$ which are counted 
with weight $(-1)^{\nu_3(\pi)}$, with $\nu_3(\pi)$ denoting the number 
of multiples of 3 in $\pi$. 

Next, $Q_{6,3}(n)$ enumerates partitions $\pi$ of $n$ whose successive 
ranks take values -1, 0, or 1. By considering partitions 
$\overset\sim \to \pi \in \r$ with $\overset\sim \to \pi= \rho(\pi)$, and 
by following the reasoning in \S8, we see that to each $\overset\sim \to \pi$ 
there corresponds $\omega_4(\pi)=2^r$ partitions $\pi$ enumerated by 
$Q_{6,3}$, where $r$ is the number of even chains of $\overset\sim \to \pi$. 
Thus
$$
Q_{6,3}(n)=\sum_{\overset\sim \to \pi\in \r, \sigma(\overset\sim \to \pi)=n}
\omega_4(\overset\sim \to \pi).\tag9.5
$$
By Theorem R' we know that
$$
A_{6,3}(n)=Q_{6,3}(n).\tag9.6
$$
So from (9.4), (9.5), and (9.6) we get 

\proclaim{Theorem 4} Let $\u$ denote the set of all (unrestricted) 
partitions. Then 
$$
\sum_{\pi\in\u, \sigma(\pi)=n}(-1)^{\nu_3(\pi)}=
\sum_{\overset\sim \to \pi\in \r, \sigma(\overset\sim \to \pi)=n}
\omega_4(\overset\sim \to \pi).
$$
\endproclaim  

\head \S10. {\un{A weighted partition theorem mod 7}}\endhead

In this and the next section we will establish three results connecting 
weighted Rogers-Ramanujan partitions and partitions into parts 
$\not\equiv 0\pm i$ ($mod$ $7$), for $i=1,2,3$. The weights in all three 
cases turn out to be products of Fibonacci numbers which are defined by 
$$
F_0=0, \quad F_1=1, \quad F_n=F_{n-1}+F_{n-2}, \quad \text{for} \quad n\ge2.
\tag10.1
$$
The Fibonacci numbers enter into the discussion naturally owing to

\proclaim{Lemma} Suppose there are $n$ boxes arranged in a certain order, 
and that each box can either be empty or filled. Then the number of ways 
in which no two consecutive boxes can both be empty is $F_{n+2}$. 
\endproclaim

{\bf{Proof}}: The Lemma is obviously true for $n=1$ (because in this case 
there are two ways and $F_3=2$) and $n=2$ (because in this case there 
are three ways and $F_4=3$). 

Let the Lemma be true for $n=1,2,...,k$. Now consider $k+1$ boxes. 

{\un{Case 1}}: Box numbered $k+1$ is non-empty.

In this case the number of ways of filling the first $k$ boxes is $F_{k+2}$. 

{\un{Case 2}}: Box numbered $k+1$ is empty.

In this case the box numbered $k$ must be non-empty, and the number of 
ways of filling the first $k-1$ boxes is $F_{k+1}$. 

So the total number of ways of filling the $k+1$ boxes is
$$   
F_{k+1}+F_{k+2}=F_{k+3}
$$
by (10.1). Hence the lemma has been proved by induction. 

Next, we define a string $\psi$ in a Rogers-Ramanujan partition $\r$ to be a 
maximal sequence of parts where the difference between consecutive parts 
is $\le 3$. Thus two strings are separated by gap $\ge 4$, and every 
Rogers-Ramanujan partition can be decomposed into strings. 

Given a string $\psi$, let $\eta(\psi)$ denote the number of gaps 
equal to 3 in $\psi$. The weight $\omega_5(\psi)$ of $\psi$ is 
defined as
$$
\omega_5(\psi)=\cases F_{\eta+3}, \quad \text{if}\quad 1\in\psi, \\
                      F_{\eta+2}, \quad \text{if}\quad 1\not\in\psi.
\endcases\tag10.2
$$
The weight of a Rogers-Ramanujan partition $\overset\sim \to \pi$ is 
defined multiplicatively as 
$$
\omega_5(\overset\sim \to \pi)=\prod_{\psi}\omega_5(\psi), \tag10.3
$$
where the product is taken over all strings $\psi$ in $\overset\sim \to \pi$. 

We are now in a position to state our first mod 7 theorem. 

\proclaim{Theorem 5} Let $A_{7,3}(n)$ denote the number of partitions of $n$ 
into parts $\not\equiv 0, \pm 3$ ($mod$ $7$). Then 
$$
A_{7,3}(n)=\sum_{\overset\sim \to \pi\in \r, \sigma(\overset\sim \to \pi)=n}
\omega_5(\overset\sim \to \pi).
$$
\endproclaim

{\bf{Proof}}: We know from Theorem R that 
$$
A_{7,3}(n)=Q_{7,3}(n),\tag10.4
$$
where $Q_{7,3}(n)$ is the number of partitions $\pi$ of $n$ such that 
the successive ranks take values -1, 0, 1, or 2. We will show that 
$$
Q_{7,3}(n)=\sum_{\overset\sim \to \pi\in \r, \sigma(\overset\sim \to \pi)=n,}
\omega_5(\overset\sim \to \pi).\tag10.5
$$
Theorem 5 will follow from (10.4) and (10.5).

Consider all partitions $\pi$ enummerated by $Q_{7,3}(n)$ and the partitions 
$\rho(\pi)=\overset\sim \to \pi$ they generate. Any even part of 
$\overset\sim \to \pi$ can occur as a hook length of a certain $\pi$ with 
rank either -1 or 1. Any odd part $>1$ of $\overset\sim \to \pi$ can occur 
as a hook length of a certain $\pi$ with rank 0 or 2, but 1 as a part can 
only occur as a hook with rank 0.  

If two parts $j$ and $j+2$ of $\overset\sim \to \pi$ differ by 2, then the 
hook representation of $j$ uniquely determines the hook representation 
of $j+2$ because both hooks must have the same rank. If two consecutive parts 
of $\overset\sim \to \pi$ differ by $\ge 4$, then the hook representation 
of one has no influence on the hook representation of the other in $\pi$. 
This explains why the multiplicative formula (10.3) is true. So we need 
only determine the weights of strings $\psi$ and for this purpose we 
concentrate on gaps in $\psi$ which are exactly 3.

If two consecutive parts of $\overset\sim \to \pi$ differ by 3, then in 
their hook representations, the only disallowed rank combinations are 
-1 and 2, or 2 and -1, depending on whether the larger part is even or odd. 
For instance, if 6 and 9 are consecutive parts of $\overset\sim \to \pi$, 
then the allowable hook representations are

\bs
\line{\hfill\epsfxsize=5truein\epsfbox{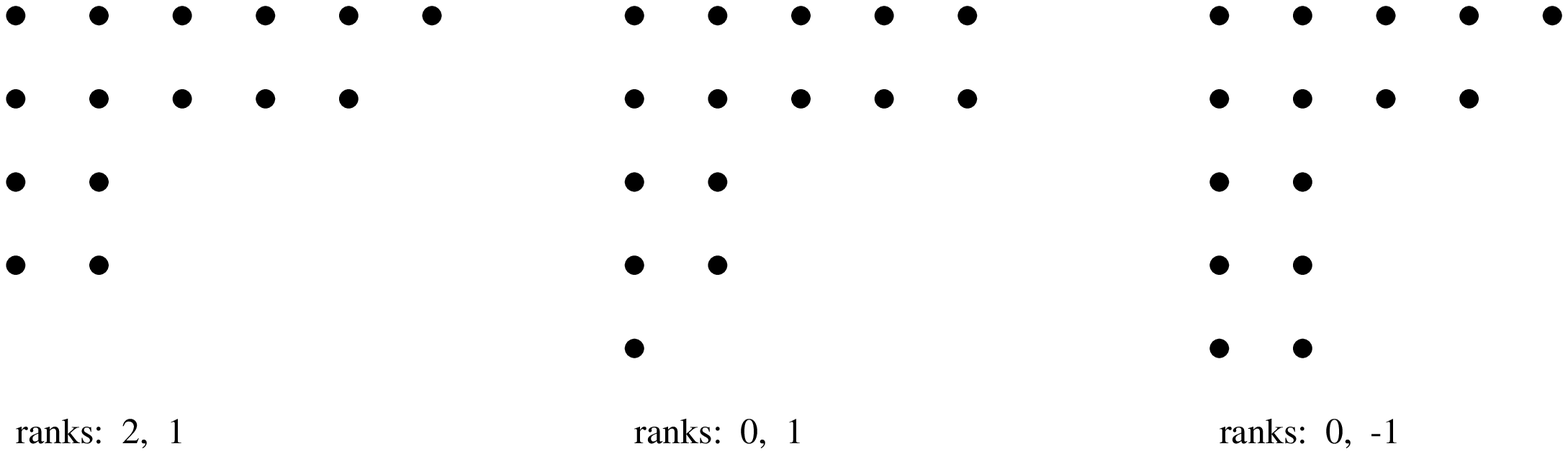}\hfill}
\centerline{\smallcaps Figure IV}
\bs

When two parts of a string differ by 3, they are of opposite parity. If two 
parts of a string differ by 2, then they are of the same parity and their 
hook representations have identical rank. Thus all hook representations 
of parts in a given chain will have the same rank once the rank of the 
hook representation of the smallest part of the chain is decided. It is the 
choice of assigning this rank to the smallest part of a chain, and 
consequently to the chain itself, that leads to weights. Thus for the purpose 
of determining the weights of strings, we may consider the decomposition 
of strings into chains. 

Consider now a string $\psi$ all of whose parts are $>1$. Let $\psi$ have 
exactly $\eta$ gaps equal to 3. This means there are $\eta+1$ chains that 
generate these $\eta$ gaps. One may think of these chains as numbered boxes. 
Adjacent chains correspond to boxes of opposite parity. If an odd chain in 
$\psi$ is represented by hooks all with rank 0, one may think of the box 
corresponding to it as being non-empty; if all hooks in the odd chain have 
rank 2, one may think of the box as being empty. Similarly, if an even chain 
is represented by hooks all with rank 1, one may think of the box as being 
non-empty, and if the hooks representing the even chain all have rank -1, 
one may think of the corresponding box as being empty. We need to assign 
ranks to the hook 
representations of these chains such that ranks 2 and -1, or -1 and 2 
cannot occur as rank combinations of adjacent chains (of hooks). This 
is the same as saying that the $\eta+1$ boxes have to be filled in such a way 
that no two adjacent boxes can be empty. From the lemma we know that 
there are $F_{\eta+3}$ ways of doing this. Thus the weight attached to 
a string $\psi$ not having 1 as a part is $F_{\eta+3}$ as in (10.2). 

Now if 1 is a part of $\psi$, then the chain in $\psi$ having 1 in it must 
have rank 0 for the hook representation of all its parts, leaving us 
no other choice. So if $\psi$ has $\eta$ gaps equal to 3, we must ignore 
the gap of 3 between the chain having 1 as a part and the next chain in 
computing $\omega_5(\psi)$. Thus we consider only $\eta-1$ gaps of 3 
and the $\eta$ chains that generate these gaps. By the reasoning of the 
preceding paragraph, the weight to be assigned to such a string will be 
$F_{\eta+2}$ as in (10.2). Thus the weights in (10.2) have been established. 

From the multiplicative definition in (10.3) it follows that to each 
$\overset\sim \to \pi\in\r$, there corresponds 
$\omega_5(\overset\sim \to \pi)$ partitions $\pi$ enumerated by $Q_{7,3}$ 
such that $\rho(\pi)=\overset\sim \to \pi$. This establishes (10.5) and so 
Theorem 5 is proved.

\head \S11. {\un{Two more weighted identities mod 7}}\endhead

The definition of weights of strings in the next two theorems will be 
a bit more complicated compared to (10.2) and the reasons will become clear 
soon. The proofs however will use methods identical to that of \S10.

Given a Rogers-Ramanujan partition $\overset\sim \to \pi$, decompose it 
into strings $\psi$ as before. Suppose $\psi$ has $\eta$ gaps equal to 3. 
Then the weight $\omega_6(\psi)$ is 
$$
\omega_6(\psi)=\cases 
F_{\eta+3}, \quad \text{if all parts of} \quad {\psi} 
\quad \text{are} \quad \ge3, \\
F_{\eta+2}, \quad \text{if either}\quad 1\in\psi \quad \text{and} 
\quad 3\not\in\psi, 
4\not\in\psi, \quad \text{or} \quad 2\in\psi, \\
F_{\eta+1}, \quad \text{if}\quad 1\in\psi \quad \text{and either} \quad 
\text{3 or 4}\in\psi.
\endcases\tag11.1
$$

As before, the weight of $\omega_6(\overset\sim \to \pi)$ is defined 
multiplicatively as
$$
\omega_6(\overset\sim \to \pi)=\prod_{\psi}\omega_6(\psi), \tag11.2
$$
We then have

\proclaim{Theorem 6} Let $A_{7,2}(n)$ denote the number of partitions of $n$ 
into parts $\not\equiv 0, \pm 2$ ($mod$ $7$). Then 
$$
A_{7,2}(n)=\sum_{\overset\sim \to \pi\in \r, \sigma(\overset\sim \to \pi)=n}
\omega_6(\overset\sim \to \pi).
$$
\endproclaim

{\bf{Proof}}: We know from Theorem R that 
$$
A_{7,2}(n)=Q_{7,2}(n).\tag11.3
$$
where $Q_{7,2}(n)$ is the number of partitions $\pi$ of $n$ such that 
the successive ranks take values 0, 1, 2, or 3. We will show that 
$$
Q_{7,2}(n)=\sum_{\overset\sim \to \pi\in \r, \sigma(\overset\sim \to \pi)=n}
\omega_6(\overset\sim \to \pi).\tag11.4
$$
Theorem 6 will follow from (11.3) and (11.4).

We consider partitions $\pi$ enumerated by $Q_{7,2}(n)$ and the partitions 
$\rho(\pi)=\overset\sim \to \pi$ they generate. If two parts of 
$\overset\sim \to \pi$ differ by $\ge 4$, then the hook representation of one 
does not influence the hook representation of the other. Thus 
$\overset\sim \to \pi$ can be decomposed into strings and 
$\omega_6(\overset\sim \to \pi)$ can be defined multiplicatively as in 
(11.2) because of this independence. 

In order to compute $\omega_6(\psi)$ for a string $\psi$, we observe that 
odd integers $>1$ can be represented by hooks with rank 0 or 2, and 
even integers $>2$ can be represented by hooks with rank 1 or 3. The 
integer 1 has a unique representation as a hook with rank 0 and 
similarly 2 has a unique representation as a hook with rank 1. If parts of 
$\psi$ differ by 3, then in their hook representation, the disallowed rank 
combinations are 0 and 3, or 3 and 0, for adjacent hooks depending on 
whether the larger part is odd or even. (The main difference between 
Theorems 6 and 5 is that in Theorem 6, the unique rank of 1 is 0, and 0 
is part of the disallowed combination, whereas in Theorem 5, the unique 
rank of 1 was 0, but the disallowed combination there did not have 1 
in it.) This is the reason for the extra complication in the definition 
of weights in (11.1).

If all parts of $\psi$ are $>4$, then clearly $\omega_6(\psi)=F_{\eta+3}$ 
as before. However, if $1\in\psi$ and either 3 or 4$\in\psi$, then 1 has 
unique rank 0, and since 0, 3, is a disallowed rank combination, the even 
chain following the chain containing 1, can only be assigned the rank 1 for 
its 
hook representation. Thus we need to discard the first two chains in computing 
weights, and so the number of chains to be considered is only $\eta-1$ 
instead of $\eta+1$. 
So by the Lemma,  the weight of $\psi$ will be $F_{\eta+1}$ as in (11.1).

If $1\in\psi$ but neither 3 nor 4 are in $\psi$, then 1 is a string by itself 
and its weight is 1. We can think of this as $F_{0+2}=F_2=1$ with $\eta=0$. 

If $2\in\psi$, then 2 admits a unique representation as a hook. This forces 
us to eliminate the chain containing 2 in computing weights and so we have 
$\eta$ chains to consider instead of $\eta+1$. So the weight in this case 
is $F_{\eta+2}$ as in (11.1). Thus we have established the weight formula 
(11.1). Therefore to each $\overset\sim \to \pi\in\r$, there corresponds 
$\omega_6(\overset\sim \to \pi)$ partitions $\pi$ enumerated by $Q_{7,2}$ 
such that $\rho(\pi)=\overset\sim \to \pi$. Thus (11.4) is established, 
and this in conjunction with (11.3) yields Theorem 6. 

To discuss the final weighted partition theorem mod 7, we consider the 
decomposition of $\overset\sim \to \pi \in\r_2$ into strings $\psi$. 
Let as before, $\eta$ denote the number of gaps equal to 3 in $\psi$. 
The weights $\omega_7(\psi)$ are defined by
$$
\omega_7(\psi)=\cases 
F_{\eta+3}, \quad \text{if all parts of} \quad {\psi} 
\quad \text{are} \quad \ge4, \\
F_{\eta+2}, \quad \text{if either 2}\in\psi \quad \text{and} \quad 
4\not\in\psi, 5\not\in\psi, \quad \text{or} \quad 3\in\psi,\\
F_{\eta+1}, \quad \text{if 2}\in\psi \quad \text{and either} \quad 
\text{4 or 5}\in\psi.
\endcases\tag11.5
$$

As always, the weight of $\omega_7(\overset\sim \to \pi)$ is defined 
multiplicatively as
$$
\omega_7(\overset\sim \to \pi)=\prod_{\psi}\omega_7(\psi), \tag11.6
$$
We then have

\proclaim{Theorem 7} Let $A_{7,1}(n)$ denote the number of partitions of $n$ 
into parts $\not\equiv 0, \pm 1$ ($mod$ $7$). Then 
$$
A_{7,1}(n)=\sum_{\overset\sim \to \pi\in \r_2, \sigma(\overset\sim \to \pi)=n}
\omega_7(\overset\sim \to \pi).
$$
\endproclaim

{\bf{Proof}}: We know from Theorem R that 
$$
A_{7,1}(n)=Q_{7,1}(n),\tag11.7
$$
where $Q_{7,1}(n)$ is the number of partitions $\pi$ of $n$ such that 
the successive ranks take values 1, 2, 3, or 4. We will show that 
$$
Q_{7,1}(n)=\sum_{\overset\sim \to \pi\in \r_2, \sigma(\overset\sim \to \pi)=n}
\omega_7(\overset\sim \to \pi).\tag11.8
$$
Theorem 7 will follow from (11.7) and (11.8).

Consider partitions $\pi$ enumerated by $Q_{7,1}(n)$ and the partitions 
$\rho(\pi)=\overset\sim \to \pi$ they generate. As before, owing to 
independence in assigning hook representations for adjacent parts of 
$\overset\sim \to \pi$ differing by $\ge 4$, we get the multiplicative 
formula (11.6). 

In order to determine the weights of strings, we observe that odd integers 
$>4$ admit hook representations with rank 2 or 4, and all even integers $>2$ 
admit hook representations with rank 1 or 3. The integer 2 has a unique 
representation as a hook with rank 1, and the integer 3 has a unique 
representation as a hook with rank 2. The integer 1 can only have a hook 
representation with rank 0, but 0 is not an allowed value of the rank. 
Thus all parts of $\overset\sim \to \pi$ are $>1$, and so 
$\overset\sim \to \pi \in\r_2$. 

If all parts of $\psi$ are $\ge4$, then we can consider all 
$\eta+1$ chains that generate the $\eta$ gaps, and the weight is 
$F_{\eta+3}$ as in (11.5). 

The disallowed rank combinations for adjacent hooklengths differing by 
3 are 1 and 4 
or 4 and 1. So if $2\in\psi$ and either 4 or 5 $\in\psi$, then the even 
chain containing 2 has rank 1 in its hook representation, and the odd 
chain following it must have rank 2 for its hook representation. So we can 
consider only $\eta-1$ chains in computing weights, and the weight in 
this case is $F_{\eta+1}$ as in (11.5). Now if $2\in\psi$ and neither 
4 nor 5 belong to $\psi$, then 2 is a string by itself with weight 1, which is 
to be interpreted as $F_{0+2}=1$ with $\eta=0$ as in (11.5). 

Finally, if $3\in\psi$, then the chain containing 3 has rank 2 for its 
hook representation and must be discarded in computing weights. So we have 
$\eta$ chains to consider, and the weight of $\psi$ by the Lemma is 
$F_{\eta+2}$ as in (11.5). 

Therefore to each $\overset\sim \to \pi\in\r_2$, there corresponds 
$\omega_7(\overset\sim \to \pi)$ partitions $\pi$ enumerated by $Q_{7,1}$ 
such that $\rho(\pi)=\overset\sim \to \pi$. Thus (11.8) is established, 
and this in conjunction with (11.7) yields Theorem 7. 
         
{\un{Remarks}}: In [3], weights which are products of Fibonacci numbers are 
attached to partitions into parts differing by $\ge4$ and these led to 
Rogers-Ramanujan partitions. Here we are attaching such weights to the 
Rogers-Ramanujan partitions and showing that these lead to partitions 
into parts $\not\equiv 0,\pm i$ ($mod$ $7$), for $i=1,2,3$. The first 
time partitions into parts $\not\equiv 0, \pm i$ $(mod$ $7)$ were discussed 
in the context of extensions of the Rogers-Ramanujan partition theorems 
was by Gordon [16]. In this paper only the congruential side mod 7 in 
Gordon's theorems are considered and not his difference conditions. 
  
\head \S12. {\un{Prospects}}\endhead

Recently we have obtained a bounded version of the G\"ollnitz partition 
theorem (see Theorem 1 of [6]), that is, a stronger form of Theorem C with 
bounds on the parts enumerated by the functions $\vee(n)$ and $C(n)$. This 
result is deduced as a consequence of a new finite identity which reduces to 
(5.1) when the bounds tend to infinity. By applying the method of \S6 
to this bounded version of the G\"ollnitz theorem, we are able to obtain the 
following new finite versions of the Jacobi triple product identity 
$$
\sum_{\ell=0}^{L}(-1)^{L+\ell}q^{2(T_L-T_\ell)}\sum_{n=-\ell}^{\ell}A^nq^{n^2}=
$$
$$
\sum_{i,j,k}(-1)^kA^{i-j}q^{2T_i+2T_j+2T_k-i-j}
\left[\matrix L-k\\i\endmatrix\right]_{q^2}
\left[\matrix L-i\\j\endmatrix\right]_{q^2}
\left[\matrix L-j\\k\endmatrix\right]_{q^2},\tag12.1
$$
and Lebesgue's identity
$$
\sum_{r,s}q^{2(T_r+T_s)}\left[\matrix L-s\\r\endmatrix\right]_{q^2}
\left[\matrix r+1\\s\endmatrix\right]_{q^2}(-1)^sA^{2s}=
$$
$$
\sum_{i,j,k}(-1)^j q^{2(T_i+T_j+T_k)-i-j}A^{i+j}
\left[\matrix L+1-k\\i\endmatrix\right]_{q^2}
\left[\matrix L+1-i\\j\endmatrix\right]_{q^2}
\left[\matrix L+1-j\\k\endmatrix\right]_{q^2}-
$$
$$
q^{2(L+1)}\sum_{i,j,k}(-1)^j q^{2(T_i+T_j+T_k)-i-j}A^{i+j}
\left[\matrix L-k\\i\endmatrix\right]_{q^2}
\left[\matrix L-i\\j\endmatrix\right]_{q^2}
\left[\matrix L-j\\k\endmatrix\right]_{q^2}.\tag12.2
$$
In (12.1) and (12.2) the symbols $\left[\matrix n\\m\endmatrix\right]_q$ are 
defined by
$$ 
\left[\matrix n\\m\endmatrix\right]_q=\frac{(q)_n}{(q)_m(q)_{n-m}}
$$
for integers $n\ge m\ge 0$. When $L$ tends to infinity, (12.1) reduces to 
(2.7), and (12.2) essentially reduces to Lebesgue's identity dilated by a 
factor of 2 (see Andrews [9], Ch. 2)
$$
(1-Aq^2)\sum_r\frac{q^{2T_r}(A^2q^4;q^2)_r}{(q^2;q^2)_r}=
(-q^2;q^2)_\infty(A^2q^2;q^4)_\infty=
\frac{(A^2q^2;q^4)_\infty}{(q^2;q^4)_\infty}\tag12.3
$$ 
from which Theorem S follows. In deriving (12.1) and (12.2), special 
attention must be paid to the ordering (6.3). This is because the bounds 
on the parts in Theorem A lead to certain exceptional cases at the boundary 
when the transformations (6.2) are applied (see [7] for details).

In the second part of this paper we have concentrated on hooks that have 
vertices on the main diagonal of the Durfee square, and hook differences 
(successive ranks) that 
take either 3 or 4 consecutive integer values. When the number of values 
taken by the rank is 3, the weights turned out to be powers of 2 as in 
Sections 8 and 9. When the successive ranks took four consecutive integer 
values as in Sections 10 and 11, the weights were products of Fibonacci 
numbers. One way to generalize this is to consider successive ranks taking
 more integer values and discuss the weighted identities they lead to. 
With this in mind, we have recently investigated the case of 5 successive 
integer values for the ranks, and even here the weights are more intricate 
than the ones considered in this paper. In this situation we need to 
decompose partitions in $\r$ into blocks of parts differing by $\le4$. 
If such a block has all gaps $\le 3$, then its weight is of the form a power 
of 2 times a power of 3, but otherwise the determination 
of the weights is more involved. 

An even more challenging question is to consider partitions with prescribed 
hook differences where the hook vertices are not on the main diagonal as in 
[10], and discuss the weighted partition theorems they lead to. We plan to 
consider such questions in the future.

{\bf{Acknowledgements}}: We would like to thank Frank Garvan for help with 
the diagrams.

\enddocument